\documentclass[]{siamart171218}

\usepackage{lipsum}
\usepackage{amsfonts}
\usepackage{graphicx}
\usepackage{epstopdf}
\usepackage{algorithmic}
\ifpdf
  \DeclareGraphicsExtensions{.eps,.pdf,.png,.jpg}
\else
  \DeclareGraphicsExtensions{.eps}
\fi

\newsiamremark{remark}{Remark}
\newsiamremark{hypothesis}{Hypothesis}
\crefname{hypothesis}{Hypothesis}{Hypotheses}
\newsiamthm{claim}{Claim}

\headers{A Runge-Kutta Method of Order Six for Lawson Integration}{M. Golden}

\title{An Explicit Sixth Order Runge-Kutta Method for Simple Lawson Integration\thanks{Submitted to the editors December 1, 2025.}}

\author{Matthew Golden\thanks{School of Physics, Georgia Institute of Technology, Atlanta, GA}, \email{matthew.golden@gatech.edu}}
\usepackage{amsopn}

\makeatletter
\newcommand*{\addFileDependency}[1]{
  \typeout{(#1)}
  \@addtofilelist{#1}
  \IfFileExists{#1}{}{\typeout{No file #1.}}
}
\makeatother

\newcommand*{\myexternaldocument}[1]{%
    \externaldocument{#1}%
    \addFileDependency{#1.tex}%
    \addFileDependency{#1.aux}%
}

\ifpdf
\hypersetup{
  pdftitle={A Runge-Kutta Method of Order Six for Lawson Integration},
  pdfauthor={M. Golden}
}
\fi

\myexternaldocument{ex_supplement}

\usepackage{subcaption}
\usepackage{tikz}
\usepackage{fancyvrb}
\usepackage{xcolor}

\definecolor{MATLABcomment}{RGB}{2, 128, 9}
\definecolor{MATLABkeyword}{RGB}{14, 0, 255}
\begin{document}

\maketitle

\begin{abstract}
  Explicit Runge-Kutta schemes become impractical when a stiff linear operator is present in the dynamics. This failure mode is quite common in numerical simulations of fluids and plasmas. Lawson proposed Generalized Runge-Kutta Processes for stiff problems in 1967, in which the stiff linear operator is treated fully implicitly via matrix exponentiation. Any Runge-Kutta scheme induces valid Lawson integration, but a scheme is exceptionally simple to implement if the abscissa $c_i$ are \textit{ordered and equally spaced}. Classical RK4 satisfies this requirement, but it is difficult to derive efficient higher order schemes with this constraint. Here I present an explicit sixth order method identified with Newton-Raphson iteration that provides simple Lawson integration.
\end{abstract}



\section{Introduction}
Explicit Runge-Kutta methods have enabled the numerical exploration of high dimensional systems of ordinary differential equations (ODEs). Such high dimensional dynamics commonly arise from the spatial discretization of partial differential equations (PDEs) in applications of fluid dynamics and plasma physics. These applications commonly lead to autonomous\footnote{Non-autonomous dynamics in $\mathbb{R}^n$ are autonomous dynamics in $\mathbb{R}^{n+1}$.} ODEs of the form
\begin{equation}
    \dot{\bf u} = {\bf g}({\bf u}) + \hat{A} {\bf u},
\label{eq:dynamics}
\end{equation}
where ${\bf u}\in\mathbb{R}^n$ is the state vector, ${\bf g}: \mathbb{R}^n \rightarrow \mathbb{R}^n$ is a smooth nonlinear function, and $\hat{A} \in \mathbb{R}^{n\times n}$ is a linear operator. The operator $\hat{A}$ usually corresponds in some part to dissipation, so it is expected to have eigenvalues $\lambda$ with $\textrm{Re}(\lambda) \ll 0$. If the timescales of these decaying modes are substantially larger than the timescales of the nonlinear term ${\bf g}$, then the problem is \textit{stiff}. 

There are a wide variety of numerical methods for stiff problems of the form \eqref{eq:dynamics}. Perhaps the most accessible are the implicit-explicit (IMEX) Runge-Kutta schemes, which use two sets of coefficients for the Runge-Kutta process \cite{ascher1997implicit}. The first set is explicit and for ${\bf g}$, while the second set is usually diagonally implicit and used to handle $\hat{A}$. IMEX schemes require repeatedly solving the system of equations $(\hat{I} - a_{ii}h \hat{A})$. Another approach is the use of Exponential Runge-Kutta (ExpRK), which uses the integral identity
\begin{equation}
    {\bf u}(t) = e^{t \hat{A}} {\bf u}_0 + \int_0^t ds \, e^{(t-s)\hat{A}} {\bf g}({\bf u}(s))
    \label{eq:dynamics2}
\end{equation}
to replace the standard Runge-Kutta coefficients $a_{ij}$ and $b_i$ with a collection of analytic functions $a_{ij}(z)$ and $b_i(z)$ that are evaluated for various $z \propto h \hat{A}$ \cite{hochbruck2005explicit, hochbruck2005exponential}. Such schemes include  ETDRK4 of Cox \& Matthews \cite{cox2002exponential}. While these ExpRK schemes have nice topological guarantees like the preservation of equilibria, their implementation requires careful evaluation of many analytic functions. 

Lawson investigated such stiff systems in 1967 \cite{lawson1967generalized} and was inspired to apply explicit Runge-Kutta methods directly to the integral of  Equation \eqref{eq:dynamics2}. Lawson called such a scheme a \textit{Generalized Runge-Kutta Process} and presented the algorithm
\begin{align}
    & {\bf k}_i = h{\bf g}({\bf U}_i), \nonumber \\
    & {\bf U}_i = e^{c_i h \hat{A}} {\bf u}_0 + \sum_j a_{ij} e^{(c_i - c_j)h \hat{A}} {\bf k}_j,\nonumber \\
    & {\bf u}^+ = e^{h\hat{A}}{\bf u}_0 + \sum_i b_i e^{(1-c_i) h \hat{A}} {\bf k}_i. \label{eq:GenRK}
\end{align}
Lawson Runge-Kutta schemes have been utilized for the numerical solution of the linear and nonlinear  Schr{\"o}dinger equation \cite{berland2005solving, cano2015projected, besse2017high, hochbruck2020convergence}, stochastic differential equations \cite{debrabant2021runge}, reaction-diffusion-reaction equations \cite{caliari2024accelerating}, Vlasov equations \cite{boutin2024modified, crouseilles2024exponential}, the nonlinear Maxwell equations \cite{pototschnig2009time}, and 3D plasma dynamics \cite{wang2020exponential}.

For an arbitrary explicit Runge-Kutta scheme with $s$ stages, applying such a scheme requires the computation of $O(s)$ matrix exponentials of the form $\exp( (c_{i+1} - c_i) h \hat{A} )$. These matrix exponentials can be computed efficiently if one is in the eigenbasis of $\hat{A}$ or approximated with Krylov subspace methods \cite{hochbruck1997krylov}. The management of these exponentials is a major barrier to the implementation of Lawson integration for an arbitrary Runge-Kutta scheme. Lawson noted from experiments with RK4 that ``\textit{if distinct $c_i$ are equally spaced in $[0,1]$ and $0 = c_1 \leq c_2 \leq \cdots \leq c_s \leq 1,$ the formulas are computationally much simpler}" \cite{lawson1967generalized}. Lawson never explicitly states how this constraint on $c_i$ simplifies the implementation of Equations \eqref{eq:GenRK}. In fact, Lawson immediately obfuscated his apparent insight by writing in the next paragraph ``\textit{We also need an algorithm for the computation of matrices $\exp(c_ihA)$.}" \cite{lawson1967generalized}. 

What Lawson must have noticed for RK4 (but failed to explicitly state) is that if $c_{i+1} = c_i +\Delta c$ (constant increment) or $c_i$ (no increment) for all $i$, then Equations \eqref{eq:GenRK} can be implemented with the evaluation of a single matrix exponential $\exp(\Delta c h \hat{A})$. I will distinguish this case from general Lawson integration with the name \textit{Simple Lawson Runge-Kutta} (SLRK) as described in Algorithm \ref{alg:LRK}.

\begin{algorithm}
\caption{Simple Lawson Runge-Kutta}
\label{alg:LRK}
\begin{algorithmic}[1]
\STATE \textbf{Given:} initial condition $\mathbf{u}$, timestep $h$, and number of timesteps $n$
\STATE $\hat{e} \gets \exp\!\bigl(\Delta c \, h \hat{A}\bigr)$
\FOR{$i=1$ \textbf{to} $n$}
  \STATE $\mathbf{k}_1 \gets h\, \mathbf{g}(\mathbf{u})$
  \FOR{$j=2$ \textbf{to} $s$}
    \IF{$c_{j} - c_{j-1} = \Delta c$}
      \STATE $\mathbf{u} \gets \hat{e}\,\mathbf{u}$
      \STATE ${\bf k}_m \gets \hat{e} {\bf k}_m$ {\bf for all} $m < j.$
    \ENDIF
    \STATE $\mathbf{k}_j \gets h\, \mathbf{g}( {\bf u} + \sum_j a_{ij} {\bf k}_j )$
  \ENDFOR
  \STATE ${\bf u} \gets {\bf u} + \sum_i b_i {\bf k}_i$
\ENDFOR
\STATE \RETURN $\mathbf{u}$
\end{algorithmic}
\end{algorithm}

To make this algorithm more intuitive, I will present a MATLAB implementation of SLRK4 in which some of the control flow of Algorithm \ref{alg:LRK} is unrolled. The Algorithm \ref{alg:SLRK4} is precisely standard RK4 interleaved with multiplication by $\exp(\Delta c h\hat{A})$ every time $c_i$ is incremented by $\Delta c$.
\begin{algorithm}
\caption{Simple Lawson Runge-Kutta 4 (SLRK4)}
\label{alg:SLRK4}
\begin{Verbatim}[commandchars=\\\{\}]
{\color{MATLABkeyword}function} u = SLRK4(u,t,n,g,A)
  {\color{MATLABcomment}%Numerically solve du/dt = g(u) + Au}
  {\color{MATLABcomment}%u - initial condition}
  {\color{MATLABcomment}%t - final time}
  {\color{MATLABcomment}%n - number of timesteps}
  {\color{MATLABcomment}%g - function handle}
  {\color{MATLABcomment}%A - matrix}
  
  h = t/n; {\color{MATLABcomment}%timestep}
  e = expm(h/2 * A); {\color{MATLABcomment}%The lone exponential we will need}
  {\color{MATLABkeyword}for} i = 1:n
    k1 = h*g(u);
    
    {\color{MATLABcomment}%c increments 0->1/2}
    u=e*u; k1=e*k1;
    
    k2 = h*g(u + k1/2);
    k3 = h*g(u + k2/2);
    
    {\color{MATLABcomment}%c increments 1/2->1}
    u=e*u; k1=e*k1; k2=e*k2; k3=e*k3; 
    
    k4 = h*g(u + k3);
    u = u + k1/6 + k2/3 + k3/3 + k4/6;
  {\color{MATLABkeyword}end}
{\color{MATLABkeyword}end}
\end{Verbatim}
\end{algorithm}

The requirement of equally spaced and ordered $c_i$ only becomes nontrivial for schemes with three or more stages. The Heun method of order three \cite{heun1900neue} has $\Delta c = 1/3$, and the classical RK4 has $\Delta c = 1/2$. Lawson used RK4 to generate a his scheme \cite{lawson1967generalized}. A year prior, Lawson published a family of six-stage fifth order methods \cite{lawson1966order} that can satisfy the constraints on $c_i$ with $\Delta c = 1/4$. While many explicit sixth order Runge-Kutta schemes have been published by 
Huta (1956) \cite{Huta1956}, 
Butcher (1964) \cite{butcher1964runge}, 
Lawson (1967) \cite{lawson1967order},
Luther (1968) \cite{luther1968explicit}, 
Fehlberg (1968) \cite{fehlberg1968classical},
Prince \& Dormand (1981) \cite{prince1981high},
Verner (1991) \cite{verner1991some},
Papakostas et al. (1996) \cite{papakostas1996general}, and others, none appear to satisfy the requirement that $c_i$ be ordered and equally spaced to allow simple Lawson integration. The purpose of this manuscript is to provide such a scheme. 

\section{Main Results}

A sixth order Runge-Kutta scheme must satisfy 37 order conditions \cite{butcher1963coefficients}. Butcher proved by construction that explicit solutions to these equations exist with seven stages. Constructing such a solution by hand is only feasible by assuming clever redundancy in the order conditions. These redundancy conditions are not aligned with the desire for equally spaced and ordered $c_i$; it is unlikely that convenient high order schemes with this constraint will be found by hand. Instead, I turn to Newton-Raphson iteration
\begin{equation}
    {\bf x}_{n+1} = {\bf x}_n - \gamma{\bf J}_n^{+} {\bf F}_n,
\end{equation}
where ${\bf x} = b_i \oplus a_{ij}$ is the set of Runge-Kutta coefficients, ${\bf F}_n = {\bf F}({\bf x}_n)$ is the 37 order conditions and the deviation of $c_i$ from a prespecified evenly spaced grid, ${\bf J}^+$ is the pseudoinverse of the Jacobian of $\bf F$, and $\gamma\in(0,1]$ is a damping parameter. Newton-Raphson iteration was applied to tens of thousands of random initial guesses for a seven-stage methods with $\Delta c=1/6$ and $\Delta c = 1/5$. No initial guesses converged to a root. This is strong numerical evidence that such a seven stage method does not exist with our constraints on $c_i$. However, Newton-Raphson iteration of eight-stage methods rapidly produced many viable roots of ${\bf F}$. A set of numerical coefficients with small magnitude was selected and numerically fine-tuned to make some coefficients zero or unity. The remaining numerical coefficients appeared to be simple rational numbers of small magnitude. The resulting Butcher table is Table \ref{tab:sixth_order}.
\begin{table}[h]
    \centering
    \begin{tabular}{c|cccccccc}
    0   &  \\
    1/6 &  1/6 \\
    1/6 &  1/12 & 1/12\\
    2/6 &  0 & -4/33 & 5/11\\
    3/6 &  -1/4 & -29/44 & 31/22 \\
    4/6 &  3/11 & 8/33 & -4/11 & 1/11 & 14/33\\
    5/6 &  -17/48 & -5/12 & 1 & 1 &-13/12&11/16\\
    1   &  20/39 & 12/39 & -31/39 & -1/39 & 34/39 & -11/39 & 16/39\\
    \hline
    & 13/200 & 0 & 4/25 & 11/40 &0 & 11/40 & 4/25 & 13/200
    \end{tabular}
    \caption{The Butcher table for an eight-stage, sixth order Runge-Kutta scheme. This scheme supports simple Lawson integration as the $c_i$ are ordered and equally spaced with $\Delta c= 1/6.$}
    \label{tab:sixth_order}
\end{table}
These coefficients were identified with double precision arithmetic and mapped by eye to rational numbers. It was subsequently verified by a symbolic calculation in MATLAB that this rational Butcher table exactly solves all 37 order conditions.

As fluid dynamics was the motivating problem for this scheme, I verify convergence by numerically solving the 2D incompressible Navier-Stokes equations
\begin{align}
    & \nabla \cdot {\bf u} = 0, \nonumber \\
    & \partial_t {\bf u} = -{\bf u} \cdot \nabla {\bf u} - \nabla p + \nu \nabla^2 {\bf u} + {\bf f},
    \label{eq:NavierStokes}
\end{align}
where ${\bf u}$ is the flow velocity, $p$ is the pressure, $\nu$ is the kinematic viscosity, and ${\bf f}$ is a body forcing term. I solve this PDE on a periodic domain $x,y\in[0, 2\pi]$ with $\nu = 10^{-2}$ and ${\bf f} = \sin 4 y \hat{i}$ as was studied in \cite{chandler2013invariant}. The mean flow $\langle {\bf u} \rangle_{xy}$ is an integral of motion, and it is taken to be zero. I over resolve the flow on a $1024\times1024$ periodic lattice with standard spectral methods. In practice, Equations \eqref{eq:NavierStokes} are easier to solve via the evolution of the vorticity $\omega \equiv \nabla \times {\bf u}$. The details of the vorticity-streamfunction formulation and the corresponding spectral discretization are beyond the scope of this manuscript, but can be found in \cite{peyret2002spectral}. The linear operator is taken to be $\hat{A} = \nu \nabla^2$. The spectrum of the operator is clear when $\omega$ is written as a Fourier series: $\hat{A}$ has large negative eigenvalues $-\nu(k_x^2 + k_y^2)$.

\begin{figure}[th]
    \centering
    \begin{subfigure}{0.45\textwidth}
        \includegraphics[width=\linewidth]{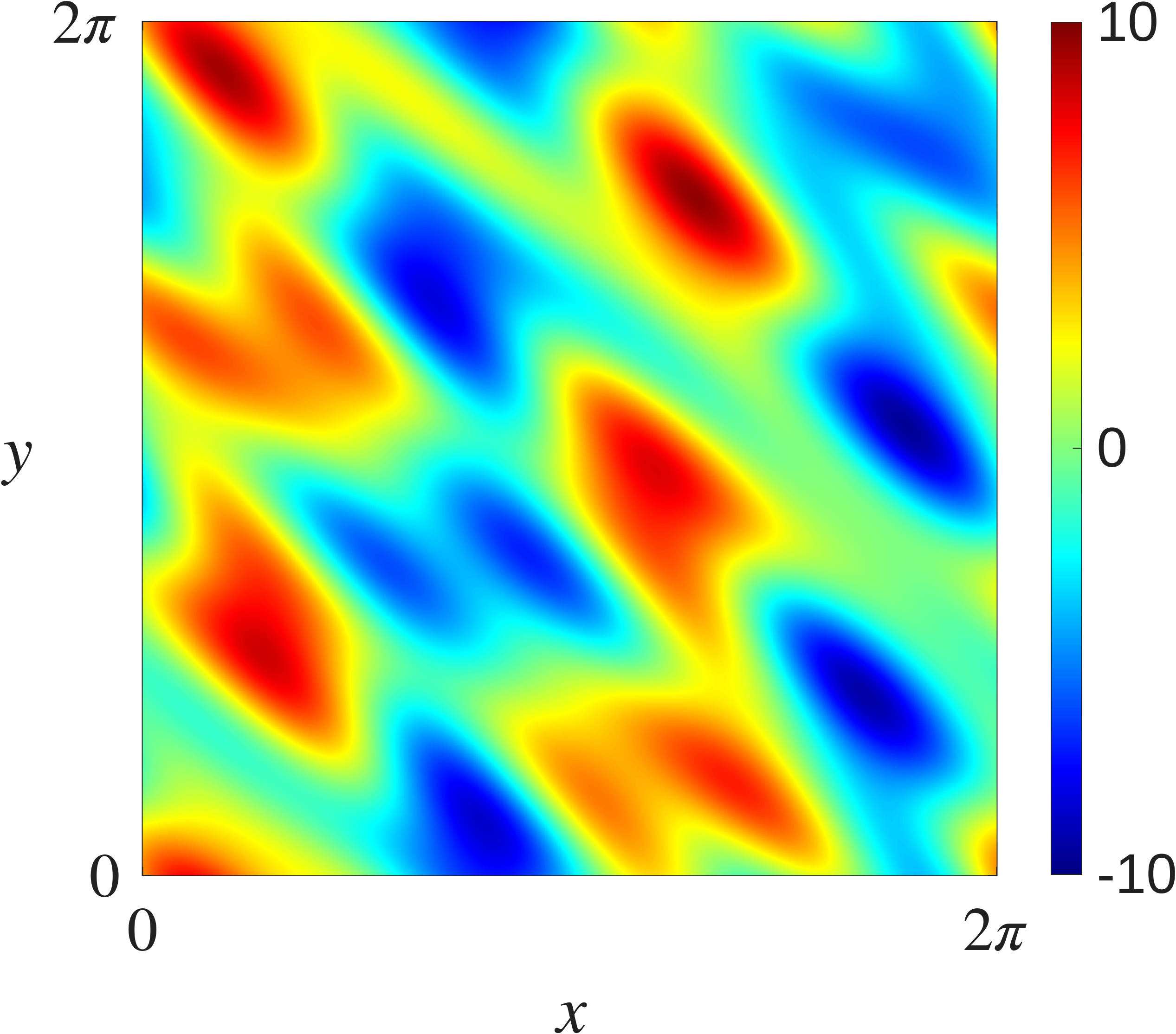}
        \caption{}
    \end{subfigure}\,
    \begin{subfigure}{0.45\textwidth}
        \includegraphics[width=\linewidth]{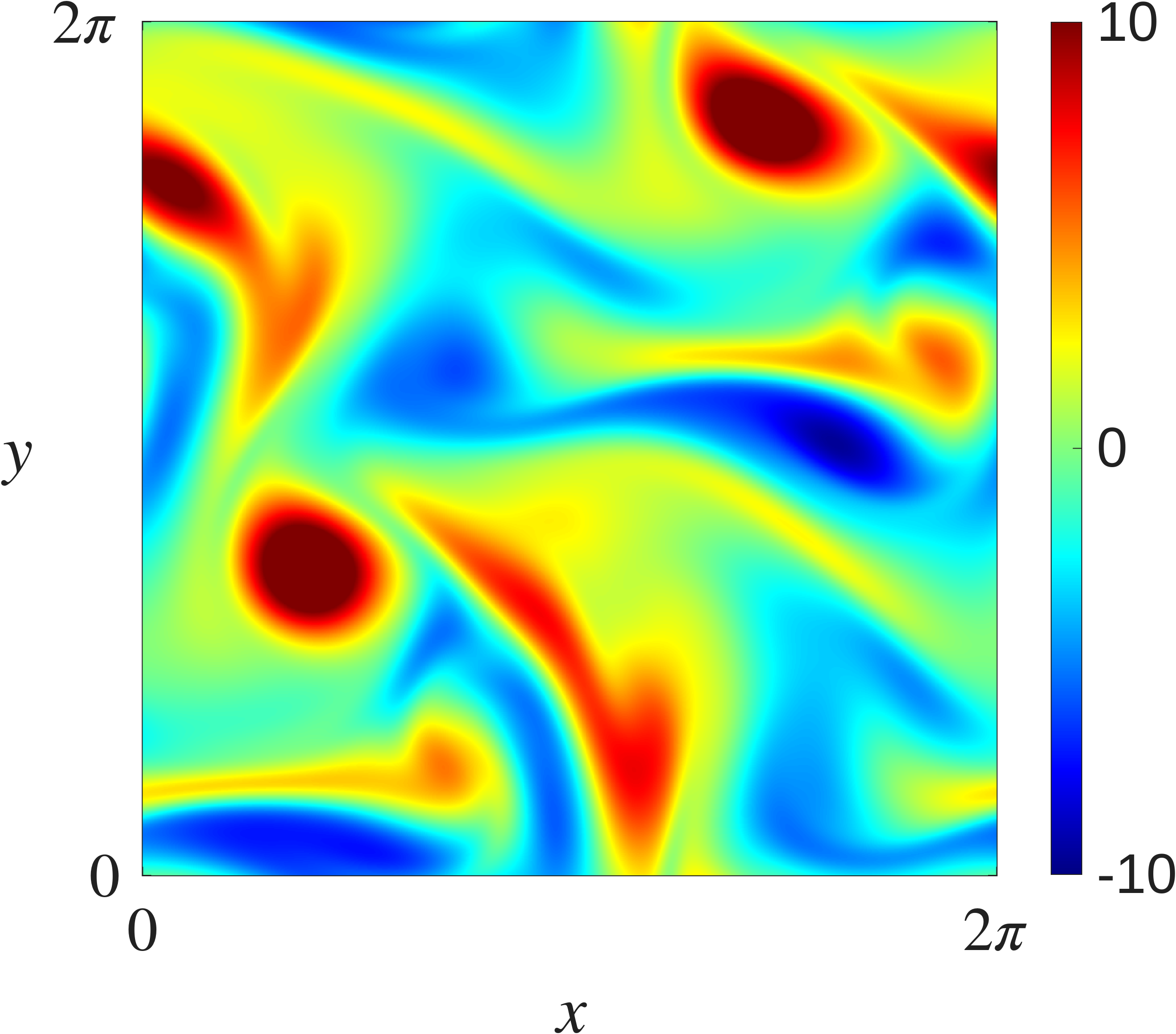}
        \caption{}
    \end{subfigure}\\[1ex]
    \begin{subfigure}{0.45\textwidth}
        \includegraphics[width=\linewidth]{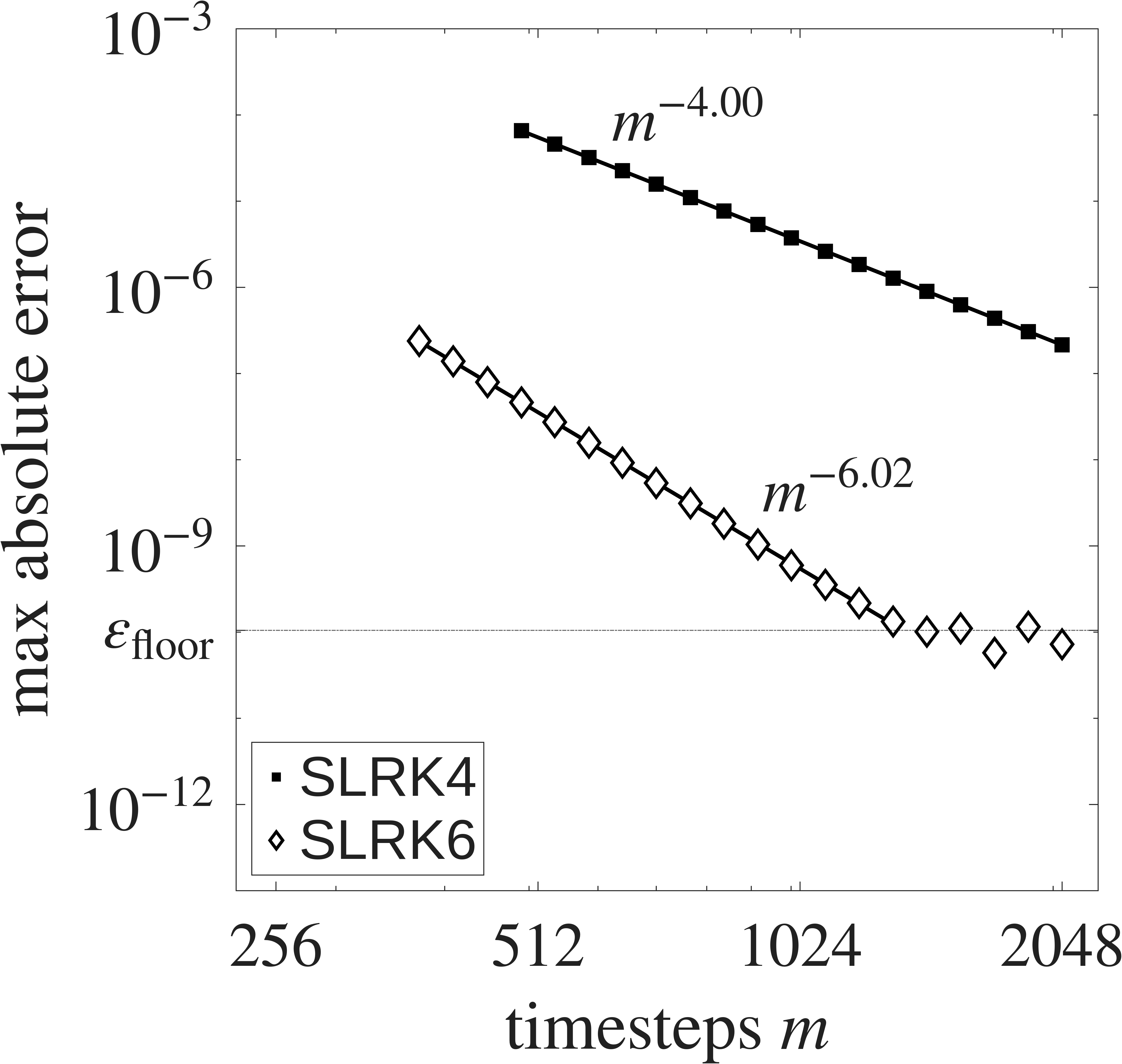}
        \caption{}
    \end{subfigure}
    \caption{(a) The analytic initial condition $\omega_0$. (b) The vorticity at $t=5$ after substantial nonlinear evolution. (c) The maximum pointwise error ($\ell_\infty$ norm) as a function of the number of timesteps $m$. A state evolved with SLRK6 and $m=2^{13}$ is treated as the ground truth. Power law fits to the error support the claimed convergence rates of SLRK. $\varepsilon_{\textrm{floor}}$ is an empirical error floor of order $10^{-10}$ likely arising from the dynamic amplification of roundoff error.}
    \label{fig:convergence_NS}
\end{figure}

An analytic initial condition\footnote{The initial state is $\omega_0 = 4\sin(2x) + 3\cos(x + 3y+0.13) + 2\sin(4x+2y+0.31) + \sin(5x + 6y + 1.23).$} without symmetry is used to test the numerical convergence of SLRK generated by the sixth order scheme of Table \ref{tab:sixth_order} and compared the scheme generated by classical RK4.  The initial condition is integrated to $t=5$ where substantial nonlinear evolution of the state has occurred. Figure \ref{fig:convergence_NS} shows the initial and final states as well as the convergence results. SLRK6 consistently outperforms SLRK4 in both error per timestep and error per ${g}$ evaluation.

\begin{figure}[th]
    \centering
        \includegraphics[width=\linewidth]{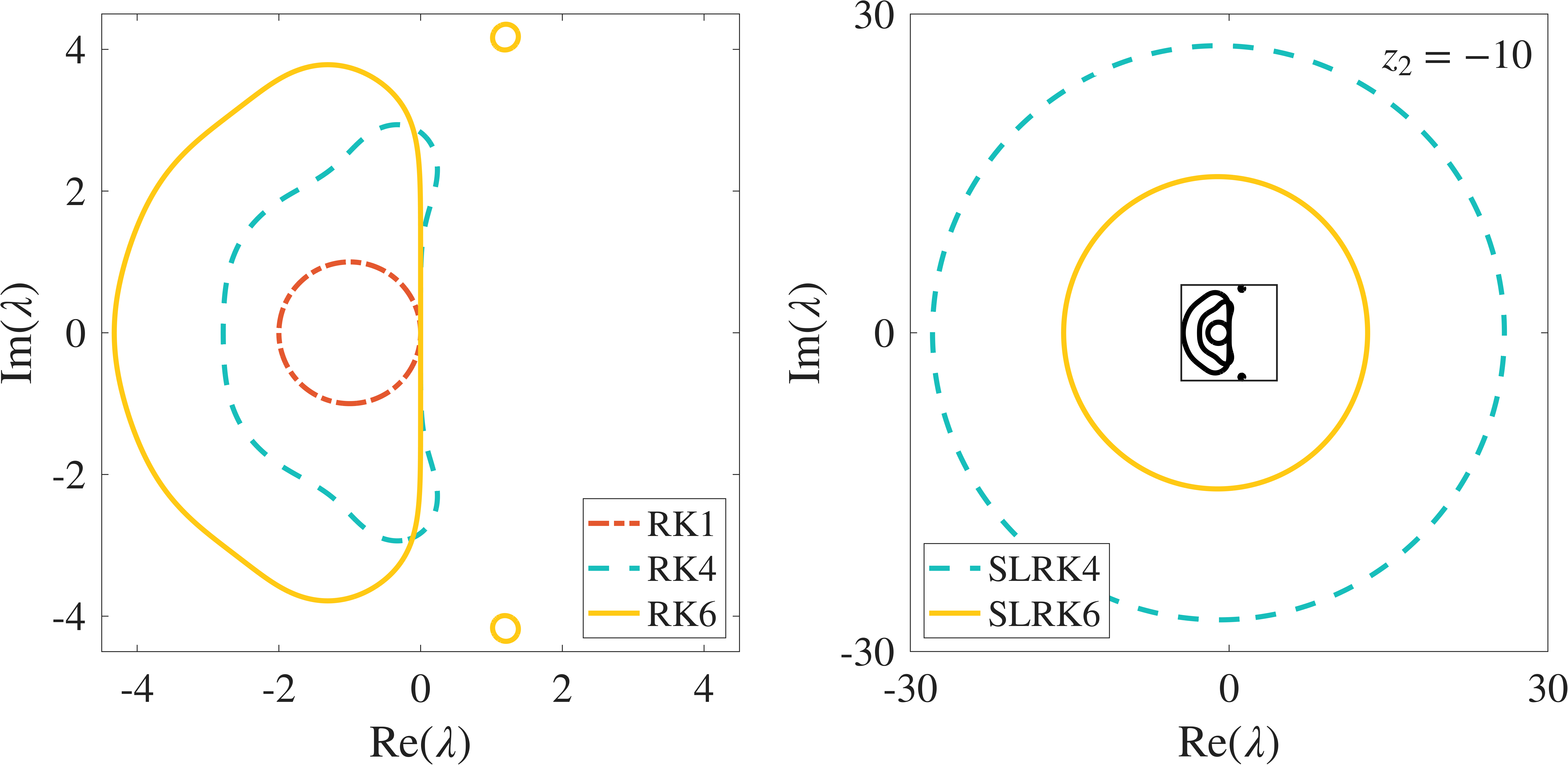}
    \makebox[0.495\linewidth]{\hfill \,\,\,\,\,\,\,(a)\hfill}
    \makebox[0.495\linewidth]{\hfill \,\,\,\,\,\,\,\,\,\,\,(b)\hfill}
    \caption{(a) Linear stability regions for RK1 (Euler's Method), RK4, and the provided RK6 scheme. (b) Linear stability regions for SLRK4 and SLRK6 with $z_2=-10$. The stability curves of subfigure (a) are shown in black to provide a sense of scale.}
    \label{fig:linear_stability}
\end{figure}

The stability function $\Phi(z)$ is a useful characterization of linear stability for Runge-Kutta methods. Numerically solution of the linear ODE $\dot{x} = \lambda x$ with timestep $h$ becomes the discrete time map $x_{n+1} = \Phi(z) x_{n}$, where $z \equiv h \lambda$. For RK4 and the RK6 method presented here, the stability polynomials are 
\begin{align}
    &\Phi_{RK4}(z) = 1 + z + \frac{z^2}{2!} + \frac{z^3}{3!} + \frac{z^4}{4!},\\
    &\Phi_{RK6}(z) = 1 + z + \frac{z^2}{2!} + \frac{z^3}{3!} + \frac{z^4}{4!} + \frac{z^5}{5!} + \frac{z^6}{6!} + \frac{29 z^7}{178200}.
\end{align}
The stability polynomials must match the power series of an exponential up to their order. These explicit schemes are numerically stable for $z$ where $|\Phi(z)| \leq 1$. These regions are sketched in Figure \ref{fig:linear_stability}. This linear stability analysis can be extended to Lawson Runge-Kutta with two timescales. Consider the ODE $\dot{x} = \lambda_1 x + \lambda_2 x$, where $\lambda_1$ is treated explicitly and $\lambda_2$ is treated implicitly with Lawson integration. Defining $z_1 \equiv \lambda_1 h$ and $z_2 \equiv \lambda_2 h$, 
\begin{align}
    &\Phi_{\textrm{SLRK4}}(z_1, z_2) = e^{z_2}\left[1 + z_1 + \frac{z_1^2}{2!} + \frac{z_1^3}{3!} + \frac{z_1^4}{4!}\right],\\
    &\Phi_{\textrm{SLRK6}}(z_1, z_2) = e^{z_2}\left[ 1 + z_1 + \frac{z_1^2}{2!} + \frac{z_1^3}{3!} + \frac{z_1^4}{4!} + \frac{z_1^5}{5!} + \frac{z_1^6}{6!} + \frac{29 z_1^7}{178200}\right].
\end{align}
That is, the only difference between the linear stability of the Lawson scheme and standard Runge-Kutta is that the exponential evolution due to $\lambda_2$ is captured exactly. Stability contours for different choices of $z_2$ are also shown in Figure \ref{fig:linear_stability}. Note that a purely imaginary $z_2$ has no impact on linear stability: high frequency oscillations pose no problem to linear stability. For a stiff Re$(z_2) \ll 0$, the decaying exponential leads to massive regions of linear stability. Note in Figure \ref{fig:linear_stability}(b) that SLRK4 becomes more stable than SLRK6 for Re($z_2) \ll 0$. This is a generic feature of Lawson integration: low order Lawson integration is more stable than high order Lawson integration in the stiff regime since $\Phi$ grows as the highest power $\Phi(z) \sim e^{z_2}z_1^{p}$ present in $\Phi(z)$ as $z_1\rightarrow \pm \infty$ and $z_2 \rightarrow - \infty$. The boundary of stability for $z_1$ should be $O(\sqrt[p]{e^{-z_2}})$ where $p=4$ for SLRK4 and $p=7$ for SLRK6. This increase in $p$ is responsible for the reordering of stability in Figure \ref{fig:linear_stability}(b). Note that regardless of $p$, the stability regions of SLRK are far larger than their standard RK counterparts. 

\section{Conclusion}
\label{sec:conclusions}

Lawson integration is a powerful method that adapts explicit Runge-Kutta methods to stiff systems. Schemes with equally spaced and ordered $c_i$ allow easy implementation via the computation of a single matrix exponential. This manuscript proposes an explicit eight-stage sixth order scheme with rational coefficients. The order of the scheme is confirmed on a stiff incompressible Navier-Stokes simulation, and the 37 order conditions are checked exactly with symbolic algebra.

This scheme was identified with Newton-Raphson iteration on the order conditions and constraints on $c_i$. Given the speed of modern computers, it is likely that this technique can find schemes of higher order well-suited for simple Lawson integration. All code used for the production of this manuscript, including a Python implementation of this scheme for both Lawson and non-Lawson integration, is provided at https://github.com/mgolden30/LawsonRK6.

\newpage
\appendix
\section{Example Implementation}
The following MATLAB code may be useful for understanding the algorithm presented here. This implementation is equivalent to Algorithm \ref{alg:LRK} using the RK6 presented in this manuscript.
\\
\begin{Verbatim}[commandchars=\\\{\}]
{\color{MATLABkeyword}function} u = SLRK6(u,t,n,g,A)
  {\color{MATLABcomment}%Numerically solve du/dt = g(u) + Au}
  {\color{MATLABcomment}%u - initial condition of size [n,1]}
  {\color{MATLABcomment}%t - final time}
  {\color{MATLABcomment}%n - number of timesteps}
  {\color{MATLABcomment}%g - function handle}
  {\color{MATLABcomment}%A - matrix}
  h = t/n;
  e = expm(h/6 * A);
  k = zeros(numel(u),8);
  a=[0,     0,     0,      0,    0,    0,     0,     0;
     1/6,   0,     0,      0,    0,    0,     0,     0;
     1/12,  1/12,  0,      0,    0,    0,     0,     0;
     0,    -4/33,  5/11,   0,    0,    0,     0,     0;
    -1/4,  -29/44, 31/22,  0,    0,    0,     0,     0;
     3/11,  8/33, -4/11,  1/11, 14/33, 0,     0,     0;
    -17/48,-5/12,  1,     1,   -13/12, 11/16, 0,     0;
     20/39, 12/39,-31/39,-1/39, 34/39,-11/39, 16/39, 0];
  b = [13/200;0;4/25;11/40;0;11/40;4/25;13/200];
  c = sum(a,2);
  Delta_c = 1/6;
  {\color{MATLABkeyword}for} i = 1:n
    k(:,1) = h*g(u);
    {\color{MATLABkeyword}for} j = 2:8
      {\color{MATLABcomment}%Check if c_j has incremented by Delta_c}
      {\color{MATLABkeyword}if} abs(c(j) - c(j-1) - Delta_c) < 1e-9
        {\color{MATLABcomment}%Apply matrix exponential}
        u = e*u;
        k(:,1:j-1) = e*k(:,1:j-1);
      {\color{MATLABkeyword}end}
      k(:,j) = h*g(u + k(:,1:j-1)*a(j,1:j-1)');
    {\color{MATLABkeyword}end}
    {\color{MATLABcomment}%Quadrature}
    u = u + k*b;
  {\color{MATLABkeyword}end}
{\color{MATLABkeyword}end}

\end{Verbatim}


\bibliographystyle{siamplain}
\bibliography{references}

\begin{thebibliography}{10}

\bibitem{ascher1997implicit}
{\sc U.~M. Ascher, S.~J. Ruuth, and R.~J. Spiteri}, {\em Implicit-explicit
  runge-kutta methods for time-dependent partial differential equations},
  Applied Numerical Mathematics, 25 (1997), pp.~151--167.

\bibitem{berland2005solving}
{\sc H.~Berland and B.~Skaflestad}, {\em Solving the nonlinear schr{\"o}dinger
  equation using exponential integrators}, tech. report, 2005.

\bibitem{besse2017high}
{\sc C.~Besse, G.~Dujardin, and I.~Lacroix-Violet}, {\em High order exponential
  integrators for nonlinear schr\"{o}dinger equations with application to
  rotating bose--einstein condensates}, SIAM Journal on Numerical Analysis, 55
  (2017), pp.~1387--1411.

\bibitem{boutin2024modified}
{\sc B.~Boutin, A.~Crestetto, N.~Crouseilles, and J.~Massot}, {\em Modified
  lawson methods for vlasov equations}, SIAM Journal on Scientific Computing,
  46 (2024), pp.~A1574--A1598.

\bibitem{butcher1963coefficients}
{\sc J.~C. Butcher}, {\em Coefficients for the study of runge-kutta integration
  processes}, Journal of the Australian Mathematical Society, 3 (1963),
  pp.~185--201.

\bibitem{butcher1964runge}
{\sc J.~C. Butcher}, {\em On runge-kutta processes of high order}, Journal of
  the Australian Mathematical Society, 4 (1964), pp.~179--194.

\bibitem{caliari2024accelerating}
{\sc M.~Caliari, F.~Cassini, L.~Einkemmer, and A.~Ostermann}, {\em Accelerating
  exponential integrators to efficiently solve semilinear
  advection-diffusion-reaction equations}, SIAM Journal on Scientific
  Computing, 46 (2024), pp.~A906--A928.

\bibitem{cano2015projected}
{\sc B.~Cano and A.~Gonz{\'a}lez-Pach{\'o}n}, {\em Projected explicit lawson
  methods for the integration of schr{\"o}dinger equation}, Numerical Methods
  for Partial Differential Equations, 31 (2015), pp.~78--104.

\bibitem{chandler2013invariant}
{\sc G.~J. Chandler and R.~R. Kerswell}, {\em Invariant recurrent solutions
  embedded in a turbulent two-dimensional kolmogorov flow}, Journal of Fluid
  Mechanics, 722 (2013), pp.~554--595.

\bibitem{cox2002exponential}
{\sc S.~M. Cox and P.~C. Matthews}, {\em Exponential time differencing for
  stiff systems}, Journal of Computational Physics, 176 (2002), pp.~430--455.

\bibitem{crouseilles2024exponential}
{\sc N.~Crouseilles and X.~Hong}, {\em Exponential dg methods for vlasov
  equations}, Journal of Computational Physics, 498 (2024), p.~112682.

\bibitem{debrabant2021runge}
{\sc K.~Debrabant, A.~Kv{\ae}rn{\o}, and N.~C. Mattsson}, {\em Runge--kutta
  lawson schemes for stochastic differential equations}, BIT Numerical
  Mathematics, 61 (2021), pp.~381--409.

\bibitem{fehlberg1968classical}
{\sc E.~Fehlberg}, {\em Classical fifth-, sixth-, seventh-, and eighth-order
  Runge-Kutta formulas with stepsize control}, vol.~287, National Aeronautics
  and Space Administration, 1968.

\bibitem{heun1900neue}
{\sc K.~Heun et~al.}, {\em Neue methoden zur approximativen integration der
  differentialgleichungen einer unabh{\"a}ngigen ver{\"a}nderlichen}, Z. Math.
  Phys, 45 (1900), pp.~23--38.

\bibitem{hochbruck2020convergence}
{\sc M.~Hochbruck, J.~Leibold, and A.~Ostermann}, {\em On the convergence of
  lawson methods for semilinear stiff problems}, Numerische Mathematik, 145
  (2020), pp.~553--580.

\bibitem{hochbruck1997krylov}
{\sc M.~Hochbruck and C.~Lubich}, {\em On krylov subspace approximations to the
  matrix exponential operator}, SIAM Journal on Numerical Analysis, 34 (1997),
  pp.~1911--1925.

\bibitem{hochbruck2005explicit}
{\sc M.~Hochbruck and A.~Ostermann}, {\em Explicit exponential runge--kutta
  methods for semilinear parabolic problems}, SIAM Journal on Numerical
  Analysis, 43 (2005), pp.~1069--1090.

\bibitem{hochbruck2005exponential}
{\sc M.~Hochbruck and A.~Ostermann}, {\em Exponential runge--kutta methods for
  parabolic problems}, Applied Numerical Mathematics, 53 (2005), pp.~323--339.

\bibitem{Huta1956}
{\sc A.~Huta}, {\em Une amélioration numérique de la méthode de
  runge--kutta--nyström des équations différentielles pour la résolution du
  premier ordre}, Acta Facultatis Rerum Naturalium Universitatis Comenianae,
  Mathematica, 1 (1956), pp.~201--244.

\bibitem{lawson1966order}
{\sc J.~D. Lawson}, {\em An order five runge-kutta process with extended region
  of stability}, SIAM Journal on Numerical Analysis, 3 (1966), pp.~593--597.

\bibitem{lawson1967generalized}
{\sc J.~D. Lawson}, {\em Generalized runge-kutta processes for stable systems
  with large lipschitz constants}, SIAM Journal on Numerical Analysis, 4
  (1967), pp.~372--380.

\bibitem{lawson1967order}
{\sc J.~D. Lawson}, {\em An order six runge-kutta process with extended region
  of stability}, SIAM Journal on Numerical Analysis, 4 (1967), pp.~620--625.

\bibitem{luther1968explicit}
{\sc H.~Luther}, {\em An explicit sixth-order runge-kutta formula}, Mathematics
  of Computation, 22 (1968), pp.~434--436.

\bibitem{papakostas1996general}
{\sc S.~Papakostas, C.~Tsitouras, and G.~Papageorgiou}, {\em A general family
  of explicit runge--kutta pairs of orders 6(5)}, SIAM journal on numerical
  analysis, 33 (1996), pp.~917--936.

\bibitem{peyret2002spectral}
{\sc R.~Peyret}, {\em Spectral methods for incompressible viscous flow},
  vol.~148, Springer, 2002.

\bibitem{pototschnig2009time}
{\sc M.~Pototschnig, J.~Niegemann, L.~Tkeshelashvili, and K.~Busch}, {\em
  Time-domain simulations of the nonlinear maxwell equations using
  operator-exponential methods}, IEEE Transactions on Antennas and Propagation,
  57 (2009), pp.~475--483.

\bibitem{prince1981high}
{\sc P.~J. Prince and J.~R. Dormand}, {\em High order embedded runge-kutta
  formulae}, Journal of computational and applied mathematics, 7 (1981),
  pp.~67--75.

\bibitem{verner1991some}
{\sc J.~Verner}, {\em Some runge--kutta formula pairs}, SIAM journal on
  numerical analysis, 28 (1991), pp.~496--511.

\bibitem{wang2020exponential}
{\sc H.~Wang, L.~Xu, B.~Li, S.~Descombes, and S.~Lanteri}, {\em An
  exponential-based dgtd method for modeling 3-d plasma-surrounded hypersonic
  vehicles}, IEEE Transactions on Antennas and Propagation, 68 (2020),
  pp.~3847--3858.

\end{thebibliography}
\end{document}